\documentclass[onefignum,onetabnum]{siamart190516}

\usepackage{amsmath}
\usepackage{amssymb}
\usepackage{xcolor}
\usepackage{mathtools}
\usepackage{float}
\usepackage{hyperref,cleveref}
\usepackage{tikz}
\usepackage[most]{tcolorbox}
\usetikzlibrary{plotmarks}
\usepackage{algorithm}
\usepackage{algpseudocode}
\usepackage{multirow}
\usepackage{tikz}
\usepackage{hhline}
\usepackage{multicol}
\usepackage{wasysym}
\usepackage{enumitem}

\newcommand{\REMOVE}[1]{}

\newcommand{\codeComment}[1]{/\hskip -.05in / #1}

\newcommand{\database}{\mathcal{B}}

\definecolor{patchcolor}{RGB}{214,220,232}

\begin{document}
\setcounter{tocdepth}{1}

\title{Compression and Reduced Representation Techniques for Patch-Based Relaxation}
\author{
  Graham Harper\thanks{Sandia National Laboratories, ({gbharpe@sandia.gov})}
\and
  Ray Tuminaro\thanks{Sandia National Laboratories, ({rstumin@sandia.gov})}
}

\maketitle

\begin{abstract}

Patch-based relaxation refers to a family of methods for solving linear systems
which partitions the matrix into smaller pieces often corresponding to groups of adjacent degrees of freedom
residing within patches of the computational domain.
The two most common families of patch-based methods are block-Jacobi and Schwarz methods,
where the former typically corresponds to non-overlapping domains and the later implies some overlap.
We focus on cases where each patch consists of the degrees of freedom within a finite element method mesh cell.
Patch methods often capture complex local physics much more effectively than simpler point-smoothers such as Jacobi;
however, forming, inverting, and applying each patch can be prohibitively expensive in terms of both storage and computation time.
To this end, we propose several approaches for performing analysis on these patches and constructing a reduced representation.
The compression techniques rely on either matrix norm comparisons or unsupervised learning via a clustering approach.
We illustrate how it is frequently possible to retain/factor less than $5\%$ of all patches and still develop a method that converges
with the same number of iterations or slightly more than 
when all patches are stored/factored.
\end{abstract}

\section{Introduction}
\label{sec: introduction}

Relaxation methods play a key role within multigrid linear solvers.  Unfortunately, unsatisfactory convergence rates are often
observed when simple relaxation schemes (e.g., point Jacobi) are employed within standard (e.g., geometric) multigrid algorithms
on matrices associated with complex partial differential equations (PDEs).  This limitation is perhaps most well-known
for saddle-point matrices, which arise from constrained PDEs systems such as the incompressible Navier-Stokes equations (INS)~\cite{MBenzi_GHGolub_JLiesen_2005a}.
However, difficulties also emerge in other scenarios, even for symmetric positive definite matrices
such as those coming from edge element discretizations of $ \nabla \times \nabla u + \alpha u $ where
$\alpha$ is a small positive constant.  For this reason, many researchers have considered sophisticated relaxation methods based
on the notion of patches, c.f.~\cite{vstar,SPVanka_1986a,ArFaWi00,he2019local,farrell2021local,SPMacLachlan_CWOosterlee_2011a,benson,scottpaper}.
Additionally, software specifically devoted to providing patch relaxation methods can be found in ~\cite{farrell2019c}.
Patch methods can be viewed as a form of overlapping Schwarz domain decomposition, though in the multigrid context the
patches or subdomains are small relative to subdomains commonly used in domain decomposition.
Specifically, patch relaxation centers on solving many small linear systems that each correspond to a subset of equations of the large matrix. 
Perhaps one of the most well-known patch relaxation schemes is the Vanka method for incompressible flow~\cite{SPVanka_1986a},
which has also been extended to magnetohydrodynamics (MHD) systems~\cite{benson,scottpaper}.
While multigrid with patch relaxation often yields satisfactory 
convergence rates for many applications including many INS and MHD simulations,
it comes at a significant cost in both storage and computation.

The main contribution of this article is the development of a family of new
approximate patch-based relaxation methods that detect and exploit structure within a linear system.
In particular, we analyze the spectral properties of patches in order to group or cluster similar patches together,
and thus reduce the total number of patches required to solve the system.
While clustering methods have been developed for matrix-variate data
using density-based approaches~\cite{gallaugher2018finite,ferraccioli2022modal},
our clustering strategy instead focuses on the spectral approximation qualities of the solver.
This enables us to only factor/store a subset of patches in the setup phase
and reuse data (stored factorizations) in the solve phase.
In some sense, this is a dimension reduction technique for linear systems,
similar in motivation to tensor decompositions \cite{kolda2009tensor} or low-rank approximations \cite{frieze2004fast} 
as it greatly reduces the necessary data storage while increasing computational speed in our case.
However, it is an approximate method, as we allow
the reconstruction to be inexact at the cost of relaxation method accuracy.
Approximate patch-based relaxation is based on the idea that a stored factorization can be used 
to approximate the true factorizations of similar patches.
Of course, it is essential that overall convergence rates do not suffer significantly. These convergence rates 
ultimately depend on the set of patches used to create the subset of stored factorizations and how 
all patch solves are approximated using the stored factorizations. A criterion is proposed for 
measuring how well a computed factorization approximates the factorization of another patch matrix.
This criterion is then used to identify patch subsets or to cluster ptatches together. Here, the
idea is to employ a single factorization for all patches within the same cluster during the algorithm's solve phase.
This paper's intent is to show the potential behind this new family of patch relaxation methods and to illustrate
some initial approaches to determine and employ patch subsets. Future work will consider more sophisticated
compression algorithms to further reduce storage costs as well as alternative algorithms that are
more efficient during the setup phase. 

In this paper, we focus on applications where the underlying discretization is high order.
In this scenario, patch relaxation can be useful in conjunction with algebraic multigrid (AMG)
due to the fact that direct application of AMG to high-order discretizations 
is often problematic.
Instead, a common alternative is to develop a composite preconditioner that applies AMG only to a low-order
discretization in conjunction with the application of patch relaxation to the high-order discretization. 
This approach is particularly useful with advanced discretizations
(e.g., discontinuous Galerkin methods) where AMG might only
be applied to a first-order continuous Galerkin finite element method,
c.f.~\cite{NLA:NLA1816,SiTuGeScCo,BECKER2018495}.
Other related applications include matrix-free approaches,
which do not store the entire sparse discretization matrix
and are essential in some applications due to storage limitations.
Matrix-free approaches are particularly attractive for
high-order discretizations as the number of matrix nonzeros may be so large
that the discretization order would need to be reduced/limited in order to actually store the entire high-order
discretization matrix.
Unfortunately, the amount of memory required to store all
the overlapping patch factors is often comparable to 
the memory needed to store the high-order discretization matrix.
Thus, storing the patch factors defeats the low-storage advantage of 
a matrix-free approach. Further, an on-the-fly approach to patch relaxation is very costly. 
On-the-fly patch relaxation requires that a factorization
be computed each time it is needed to solve a patch sub-problem.
This means that during the solve phase an $O(p_s^3)$ Gaussian elimination algorithm
must be repeatedly employed as opposed to an $O(p_s^2)$ backsolve algorithm
that is used when a factorization is already computed.
Here, $p_s$ is the patch size or number of unknowns within a patch. 
If each patch corresponds to all unknowns within a single computational cell, 
$p_s = q (p+1)^d$ when $p^{th}$-order Lagrange polynomial 
functions are used to discretize systems with $q$ PDEs on $d$-dimensional tensor (quadrilateral/hexahedral) meshes.
It follows that $p_s$ can easily exceed $500$ in many scenarios.

While this paper focuses on exploiting structure and reducing patch storage for high-order applications,
we note that there may also be further significant gains associated 
with reusing patch factorizations when bandwidth/cache performance issues
are considered on advanced computing hardware, especially on GPUs.
While it is often not possible to simultaneously keep all patch factorizations
in some type of fast-access memory with traditional patch methods, it
should be possible to keep a subset of patch factorizations
in fast-access memory during the backsolves associated with multiple right hand sides.
Finally, we remark that while we consider patch relaxation within multigrid for high-order discretizations,
patch schemes are also useful as preconditioners (without multigrid),
especially for time dependent multiphysics problems.

This paper proceeds by formally defining a traditional patch relaxation method in Section~\ref{sec: background}.
Section~\ref{sec: metrics} describes the approximate patch relaxation idea and provides an associated minimization
function that guides the compressed representation used to approximate patch inverses. In 
Section~\ref{subsec: equality} a simple algorithm is presented to identify a subset of patch
matrices whose factorizations can be used to approximate other patch matrix inverses. 
A more sophisticated algorithm is given in Section~\ref{subsec: clustering} that leverages clustering
ideas from machine learning. Numerical results are presented in Section~\ref{sec: numerical results} illustrating
the efficacy of the compressed representation when preconditioning high-order discretization problems.
Finally, Section~\ref{sec: conclusion} provides a conclusion and discusses additional work needed to mature
the basic idea.

\section{Patch-based smoothers background}
\label{sec: background}

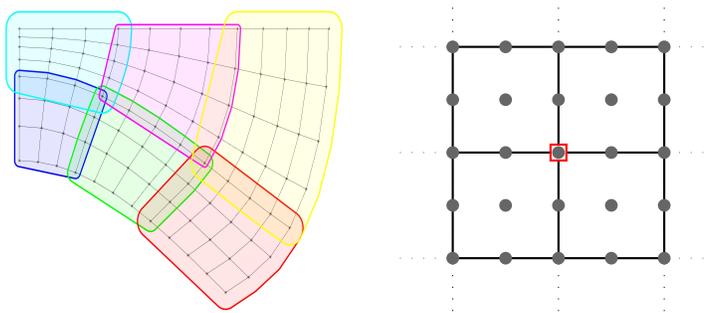
\begin{figure}[H]
\centering
\includegraphics[trim = 4.1in 4.8in 4.1in 4.0in, clip, height = 4.0cm,width = 5.5cm]{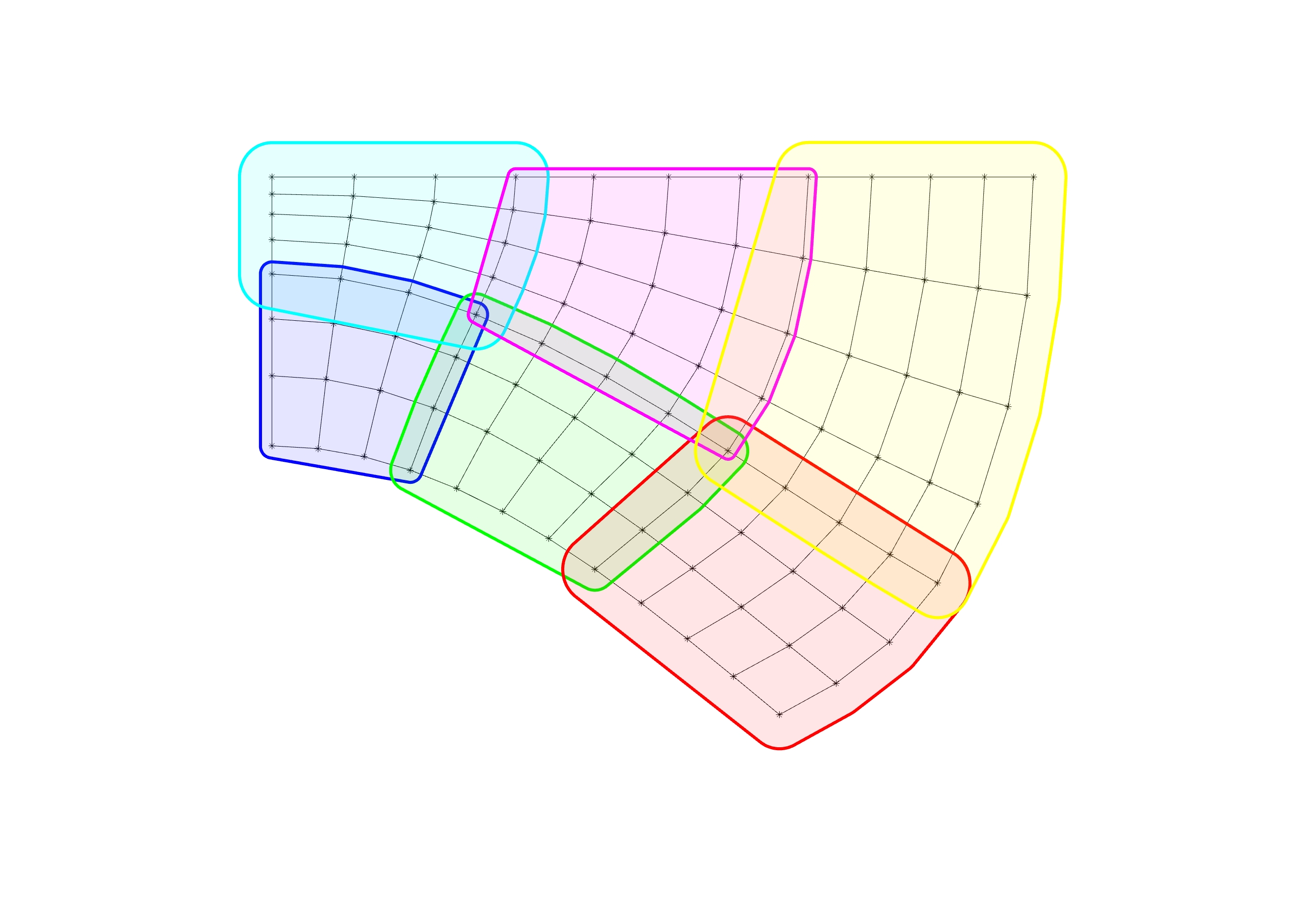}
   \newcommand{\mygrid}[3]{%
      \def\n{#1}
      \def\m{#2}

      \foreach \y in {0,\m,...,\n}{%
         \draw[gray, loosely dotted] (-1, \y) -- (\n+1, \y);
      }
      \foreach \x in {0,\m,...,\n}{%
         \draw[black, loosely dotted] (\x, -1) -- (\x, \n+1);
      }

      \foreach \y in {0,2,...,\n}{%
         \draw[black, thick] (0, \y) -- (\n, \y);
      }
      \foreach \x in {0,2,...,\n}{%
         \draw[black, thick] (\x, 0) -- (\x, \n);
      }

      \ifnum#3=1,{
         \foreach \y in {1,2,...,\n}{%
            \draw[black, densely dotted] (0, \y) -- (\n, \y);
         }
         \foreach \x in {1,2,...,\n}{%
            \draw[black, densely dotted] (\x, 0) -- (\x, \n);
         }
      }
      \else{}\fi

      \def\h{0.15}
      \draw[thick, draw=red!50!red,fill=white] (\n/2-\h, \n/2-\h) rectangle +(2*\h,2*\h);

      \foreach \y in {0,1,...,\n}{%
         \foreach \x in {0,1,...,\n}{%
            \pgfmathtruncatemacro\resultX{Mod(\x,2)==0?0:1}
            \pgfmathtruncatemacro\resultY{Mod(\y,2)==0?0:1}
            \pgfmathtruncatemacro\resultXY{\resultX*\resultY}
            \ifnum\resultX=0,{%
               \ifnum\resultY=0,{\draw[thick, draw=black!60, fill=black!60] (\x, \y) circle (2pt);}
               \else{\draw[thick, draw=black!60, fill=black!60] (\x, \y) circle (2pt);}\fi
            }
            \else{
               \ifnum\resultXY=1,{\draw[thick, draw=black!60, fill=black!60] (\x, \y) circle (2pt);}
               \else{\ifnum\resultX=1,{\draw[thick, draw=black!60, fill=black!60] (\x, \y) circle (2pt);} \else {}\fi}
               \fi
            }\fi
         }
      }
   }
   \begin{tikzpicture}[x=20pt,y=20pt]
      \begin{scope}[shift={(10pt, 10pt)}]
         \mygrid{4}{2}{0}
      \end{scope}
   \end{tikzpicture}
   \newcommand{\mycircle}  {\tikz[baseline=-0.5ex]{\draw[thick, draw=black!60, fill=black!60] (0, 0) circle (2pt);}}
\caption{Left: sample mesh with 6 patches. Right: 1 Stokes flow Vanka patch for a 
$\boldsymbol{ \mathbb{Q}}_2/\mathbb{Q}_1$
discretization that includes a pressure DoF \protect\tikz[baseline=-0.5ex]{\protect\draw[thick, draw=red!50!red,fill=white] (-0.10, -0.10) rectangle +(.2,.2);} and all velocity DoFs \protect\mycircle ~residing within elements containing the pressure DoF.  
\label{fig: sample patch} }
\end{figure}

A classical Jacobi iterative method can be viewed as a non-overlapping patch algorithm where each degree of freedom is treated as its own patch matrix.
A single Jacobi iteration applied to the $n \times n$ linear system
$$
A x = b
$$
corresponds to updating each solution unknown independently via 
\begin{align*}
  x_i^{(m+1)} 
  &=  x_i^{(m)} +
  A_{ii}^{-1} \left(b_i - \sum\limits_{j=1}^n A_{ij} x_j^{(m)} \right),
%
%
  \quad i=1,\dots,n
\end{align*}
where 
$x^{(m)}$ denotes the current approximate solution after $m$ Jacobi iterations.
When applied with damping coefficient $\omega$, Jacobi can be written succinctly as
\begin{align*}
  \mathbf{x}^{(m+1)} 
  &= 
  \mathbf{x}^{(m)} + \omega D^{-1} 
  \left(\mathbf{b} - A \mathbf{x}^{(m)} \right)
\end{align*}
where $D$ is a diagonal matrix whose only nonzeros correspond to $D_{ii} = A_{ii}$.
While the Jacobi iteration is inexpensive to store and fast to apply, it does not scale well for problems involving
complex multiphysics interactions or some high-order polynomial discretizations.
A more sophisticated block Jacobi iterative scheme is obtained by changing the definition of $D$ to instead be a block
diagonal matrix whose nonzero entries are defined by $D_{ij} = A_{ij}$. Typically, the rows/columns of 
a block might include all degrees of freedom (DoFs) that are defined at the same spatial location for a PDE system or
might include all DoFs defined at spatial locations within a small local neighborhood. 
A more general and robust patch relaxation procedure can be defined when the small neighborhoods are allowed to overlap 
as depicted on the left side of Figure~\ref{fig: sample patch}, though this can no longer be represented by a block diagonal matrix
$D$ due to the overlap. This class of algorithms is also referred to as overlapping Schwarz methods.
The right side of Figure~\ref{fig: sample patch} illustrates a single Vanka patch to solve the
two dimensional incompressible Stokes equations. Here, each patch coincides with one pressure DoF and includes all velocity
DoFs in the four elements that contain the pressure DoF.

To formally define patch relaxation, we introduce the $ p_s \times n $ boolean restriction matrix $V_k$,
$k=1,\dots,n_p$. Here, we assume that all $n_p$ patches are of the same size $p_s$ 
(often the case when preconditioning high-order discretization matrices). 
Each nonzero $V_k$ entry is given by $(V_k)_{ij} = 1$ if and only if global degree of freedom $j$ corresponds to local patch degree of freedom $i$. 
Restricting the matrix $A$ to the $k^{th}$ patch is defined by 
\begin{align*}
  A_k
  &= 
  V_k A V_k^T,
  \quad k=1,\dots,n_p.
\end{align*}
An $n \times n$  diagonal weight matrix $W = \left ( \sum\limits_{\ell=1}^{n_p} V_\ell^T V_\ell \right )^{-1}$ is defined 
to average solutions from overlap regions. This is equivalent to taking $W_{ii}$ as the reciprocal of the number of patches that include
the $i^{th}$ global DoF. 
A patch preconditioner $M^{-1}$ is 
then given by 
\begin{align}
  M^{-1}
  &=
  W \sum\limits_{k=1}^{n_p} V_k^T A_k^{-1} V_k .
  \label{eq: global patch inverse}
\end{align}
Therefore, applying $M^{-1} r$ corresponds to the restriction of $r$ to each patch, performing a patch solve, injecting the solution back to the larger problem, and then combining/averaging the results.

\textbf{Remark:} While we define patch relaxation in terms of a single patch size $p_s$, this extends to patches of different sizes. The algorithms developed later in Section~\ref{sec: metrics} may be applied to each group of same-size patches independently.

Many different patch choices are possible. As already noted, Vanka patches have proved useful for the incompressible Navier
Stokes equations. Arnold-Faulk-Winther patches~\cite{ArFaWi00} are useful for some electromagnetics applications to precondition
operators where the underlying PDE has the form $ \nabla \times \nabla u + \alpha u $. Other types of patches can be useful
for some MHD formulations.
In this paper, our definition of a patch is the smallest computational domain with exploitable regular structure.
For many discretizations and applications, this domain is exactly a mesh cell, which we recall may contain as many as
$p_s = q(p+1)^d$ degrees of freedom when solving $q$ PDEs with $p^{th}$-order Lagrange polynomials on a $d$-dimensional cube.
\begin{figure}
  \begin{tabular}{ccc}
  \resizebox{!}{3.5cm}{
  \begin{tikzpicture}
  \begin{scope}[xshift=1.5cm,yshift=0cm]
  \filldraw [fill=patchcolor, rounded corners] (0.1,0.1) rectangle (1.9,1.9);
  \draw[step=1cm,black,thin] (-1,-1) grid (2,2);
  \foreach \x in {-1,-0.8,-0.5,-0.2,0,0.2,0.5,0.8,1,1.2,1.5,1.8,2} {
    \foreach \y in {-1,-0.8,-0.5,-0.2,0,0.2,0.5,0.8,1,1.2,1.5,1.8,2} {
      \node at (\x, \y) [circle,fill,inner sep=0.7pt] {};
    }
  }
  \end{scope}
  \end{tikzpicture}}
  &
  \resizebox{!}{3.5cm}{
  \begin{tikzpicture}
  \begin{scope}[xshift=5cm,yshift=0cm]
  \filldraw [fill=patchcolor, rounded corners] (-0.3,-0.3) rectangle (1.3,1.3);
  \draw[step=1cm,black,thin] (-1,-1) grid (2,2);
  \foreach \x in {-1,-0.8,-0.5,-0.2,0,0.2,0.5,0.8,1,1.2,1.5,1.8,2} {
    \foreach \y in {-1,-0.8,-0.5,-0.2,0,0.2,0.5,0.8,1,1.2,1.5,1.8,2} {
      \node at (\x, \y) [circle,fill,inner sep=0.7pt] {};
    }
  }
  \end{scope}
  \end{tikzpicture}}
  &
  \resizebox{!}{3.5cm}{
  \begin{tikzpicture}
  \begin{scope}[xshift=5cm,yshift=0cm]
  \filldraw [fill=patchcolor, rounded corners] (-0.1,-0.1) rectangle (1.1,1.1);
  \draw[step=1cm,black,thin] (-1,-1) grid (2,2);
  \foreach \x in {-1,-0.8,-0.5,-0.2,0,0.2,0.5,0.8,1,1.2,1.5,1.8,2} {
    \foreach \y in {-1,-0.8,-0.5,-0.2,0,0.2,0.5,0.8,1,1.2,1.5,1.8,2} {
      \node at (\x, \y) [circle,fill,inner sep=0.7pt] {};
    }
  }
  \end{scope}
  \end{tikzpicture}}
  \\
  (a) Vertex-star patches
  &
  (b) Cell-centered patches
  &
  (c) Cell-restricted patches
  \end{tabular}
\caption{Comparison of different patch methods on a quadrilateral grid with a $p=4$ discretization.}
\label{fig: star vs cell patches}
\end{figure}
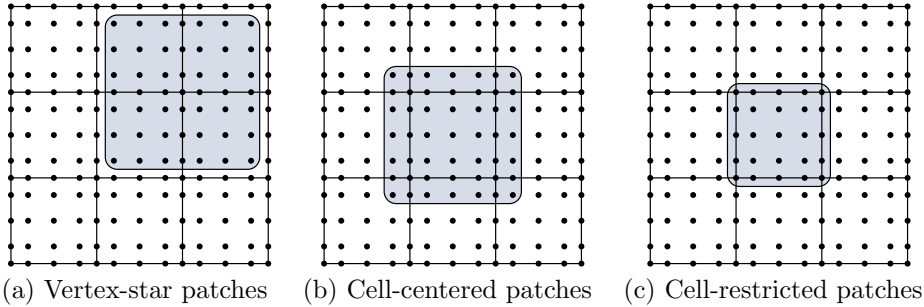

Figure~\ref{fig: star vs cell patches} highlights the different types of patches,
such as the vertex-star and cell-centered patches, which are discussed in \cite{vstar}, versus our cell-restricted patch.
Vertex-star patches are generally the most robust for patch-based relaxation methods as the polynomial degree $p$ increases;
however, a vertex-star patch may have varying sizes depending on the connectivity
of the vertex within the mesh.
This means there is irregularity in patch size on unstructured grids.
The cell-centered patch includes a halo from its neighbors,
which is generally more robust than cell-restricted patches.
We remark that stencil irregularity may be handled by grouping alike-sized
patches together before applying the algorithms in the following section,
but both vertex-star and cell-centered patches tend to also require
invasive information from the application beyond the supplied linear system.
From here forward, we utilize the cell-restricted patches shown in (c) of Figure~\ref{fig: star vs cell patches},
as we can easily detect them by searching a linear system for rows
with $q(p+1)^d$ symbolically nonzero entries for tensor-product cell shapes without requiring
additional information from the application that generated the system.
This may be adjusted accordingly depending on the cell shape and discretization.

We focus on preconditioning high-order discretizations via a combination of patch preconditioning
and multigrid applied to a first-order discretization developed on the same mesh used for the high-order discretization. Specifically,
an additive combination preconditioner has the form 
\begin{align}
M^{-1}_{combo} &= M^{-1} + P_0 M^{-1}_{amg} R_0 
  \label{eq: combination preconditioner}
\end{align}
where $M^{-1}$ is the 
patch method, $P_0$ interpolates solutions from the low order representation to the high-order representation,
$ M^{-1}_{amg}$ denotes the first order AMG preconditioner and $R_0$ restricts residual from the high-order representation to the 
low order representation.
Section~\ref{sec: numerical results} further discusses $M^{-1}_{combo}$. We close this
section by noting that there are other AMG approaches to high-order discretizations that avoid the need for a patch preconditioner.
One noteworthy possibility is to define a first order discretization using a refined mesh such that the spatial locations associated
with unknowns on both the first order and high-order systems coincide~\cite{heys2005algebraic,olson2007algebraic,10.2307/43693324}. Another possibility instead uses $p$-multigrid by creating a hierarchy of intermediate discretizations and reducing
the order of approximation by only one as we proceed through the hierarchy~\cite{Rnquist1987SpectralEM}.

\section{Approximate Patch Smoothers}
\label{sec: metrics}

To develop an approximate patch preconditioner, we need to define a subset of factored or inverse patch matrices
$$
\database = \{B_1^{-1},B_2^{-1},\dots,B_{m_p}^{-1}\} 
$$
and to identify a means of approximating each of the true patch inverses using only $\database$.
Here, $m_p (\leq n_p)$ denotes the total number of computed inverses.
While it is possible to consider sophisticated schemes that combine several $B_j^{-1}$ to approximate the $j^{th}$ patch inverse,
we restrict ourselves to the simplest case where only a single factorization is used, i.e.  $A_j^{-1} \approx B_{\phi(j)}^{-1}$.
In this case, a mapping $\phi: \{1,\dots,n_p\} \to \{1,\dots,m_p\}$ must be determined.  We will refer to $\database$ as the inverse database, 
but it may also be thought of as a compressed representation or sparse approximation of the true patch inverses $A_i^{-1}$, $i=1,\dots,n_p$.
Ultimately, $\database$ and $\phi()$ define an approximate preconditioner
\begin{align}
  \tilde{M}^{-1}
  &=
  W \sum\limits_{k=1}^{n_p} V_k^T B_{\phi(k)}^{-1} V_k  .
  \label{eq: approximate global patch inverse}
\end{align}
When $\tilde{M}$ is employed as a stand-alone preconditioner, we seek a well conditioned $\tilde{M}^{-1} A$
using a $\database$ such that $m_p$ is not too large. 
When $\tilde{M}$ is used within multigrid, we instead seek a relaxation error smoother operator $I-\omega \tilde{M}^{-1} A$ 
(with damping parameter $\omega$) that damps errors which cannot be well-represented on the first-order mesh. 
The difference between the ideal operator $I-\omega M^{-1} A$ and the approximate operator $I-\omega \tilde{M}^{-1} A$
on patch $i$ is given by
\begin{align*}
  V_i \left(\left( I-\omega M^{-1} A \right)
  -
  \left( I-\omega \tilde{M}^{-1} A \right)\right)
  V_i^T
  &=
  \omega V_i W \sum\limits_{k=1}^{n_p} \left( V_k^T \left( A_k^{-1} - B_{\phi(k)}^{-1} \right) V_k \right)AV_i^T.
\end{align*}
This results in an error term of the form $ \left(A_i^{-1} - B_{\phi(i)}^{-1}\right)A_i = I - B_{\phi(i)}^{-1} A_i $ after we have ignored leading constants and discarded terms that appear due to overlap. Therefore, a good approximation would aim to reduce $ \sum_{k=1}^{n_p} \| I - A_k B_{\phi(k)}^{-1} \|_2^2$, noting symmetric matrices may be commuted in the norm.
We consider the following minimization problem: seek a database $\mathcal{B}$ and mapping $\phi()$ minimizing $\mathcal{L}(\mathcal{B},\phi)$, where
\begin{align}
  \mathcal{L}(\mathcal{B},\phi)
  =  
  \beta |\database|
  +
  \sum_{k=1}^{n_p} 
  \| I - A_k B^{-1}_{\phi(k)} \|_2^2
  \label{eq: patch equality}
\end{align}
where the $A_k$ and $\beta > 0$ are given, $|\cdot|$ is set cardinality, and $\|\cdot\|_2$ is the matrix 2-norm.
Clearly, the second term is minimized when each $B_k^{-1}$ exactly approximates all patch inverses
assigned to it (i.e., for all $j,k$ such that $\phi(k)=j$ , $B_j^{-1} = A_k^{-1}$).
However, the first term penalizes the objective based on the size of $\database$.
Thus, the two terms in \eqref{eq: patch equality} balance the storage costs (and computation costs) associated with having a large set of approximate inverses against the quality of the inverse approximations.

\noindent
\textbf{Remark:} 
The entrywise difference Frobenius norm $\|A_k - B_{\phi(k)}\|_F^2$ or the entrywise absolute difference $\ell_1$ norm are among the cheapest to compute since these may be done before computing each $B_{\phi(k)}^{-1}$,
but using such entrywise differences also provides less information about the quality of approximation than, for example, 
comparing the spectrum of $A_k B_{\phi(k)}^{-1}$ to the identity.
For these reasons, despite the nonsymmetry, we use the 2-norm as it may be computed exactly by a small SVD or approximated by a few steps of the power method applied to $I - A_k B_{\phi(k)}^{-1}$.

\noindent
\textbf{Remark:} The minimization problem \eqref{eq: patch equality} guarantees a minimum exists, and the $\phi$ mapping must be onto to be a minimizer, but does not allow for a unique minimizer.
The minimum is guaranteed to exist as the objective function is bounded below by $\beta$ in the case of $\database$ with only one entry.
Furthermore, in the case of a minimizer, all entries of the database $\database$ must be matched with at least one $A_i$ patch,
implying the map $\phi$ must be onto.
If that were not the case, the objective function value could simply be decreased by removing unnecessary patches.
This means it is logical to consider constructive approaches for computing $\database$ in practice.
However, the minimizer is not unique because any permutation of a database that minimizes this objective function is also a minimizer.
This could be rectified by imposing additional conditions on the ordering of terms in $\database$.
These technical details are not critical here so long as the underlying minimization algorithm is not affected by permutations in the database. 

Notice that solving the above minimization problem could be computationally expensive, especially given the fact that 
$ \mathcal{L}(\mathcal{B},\phi) $ is not differentiable as $|\database|$ takes on only integer values.
However, an exact minimizer is not necessary, only a solution that gives acceptable convergence rates for a reasonable cost.
Furthermore, there are a number of ways in which an inexpensive initial guess can be used to start a minimization process. 

\subsection{Patch Equality Compression}
\label{subsec: equality}

When relatively small values of $\beta$ are considered, the $\mathcal{L}(\mathcal{B},\phi)$ minimization problem essentially attempts to exactly satisfy the first term in \eqref{eq: patch equality}. While this might lead to a large $|\database|$ for general unstructured/variable coefficient PDE 
problems, this is not the case for homogeneous coefficient PDEs on regular meshes, 
as the majority of the patch matrices (excluding boundary conditions) are identical.
This type of small $\beta$ or patch equality compression scenario is also beneficial for problems with piecewise-constant variations in materials.
These problems come up frequently when modeling geologic or crystalline structure of materials or in semiconductor simulations. One simple approach that avoids the explicit $\mathcal{L}(\mathcal{B},\phi)$ minimization problem is given in Algorithm~\ref{algo: tolerance}. This algorithm requires $O(n_p m_p) $ operations as each patch must be checked against the current $\database$ to see if there is already a suitable match or whether $\database$ should be enlarged. 
\begin{algorithm}[H]
\caption{Patch Tolerance Construction of $\database$ and $\phi()$} \label{algo: tolerance}
\begin{algorithmic}[1]
\State \textbf{Input:}   patches $\{A_1,A_2,\dots,A_{n_p}\}$, tolerance $\varepsilon$
\State $\mathcal{B}:=\{\}$, $\vec{\phi}=\vec{0}$
\For{$i=1,2,\dots,n_p$}
  \State match:=false;
  \For{$j=1,2,\dots,|\mathcal{B}|$}
    \If{$\|I-A_i B_j^{-1}\|_2<\varepsilon$}~~~~~~~~~~~\codeComment{If a suitable match is found...}
      \State match=true, $\phi(i)=j$, break;~~~\codeComment{save the index, skip to next patch}
    \EndIf
  \EndFor
  \If{match==false}~~~~~~~~~~~~~~~~~~~~\codeComment{If no suitable match is found...}
    \State append $A_i^{-1}$ to $\mathcal{B}$, $\phi(i)=|\mathcal{B}|$;~~~~~~~\codeComment{append database, $m_p$ incremented}
  \EndIf
\EndFor
\State \textbf{Output:} $\mathcal{B}=\{B_1^{-1},B_2^{-1},\dots,B_{m_p}^{-1}\}$, $\vec{\phi}$
\end{algorithmic}
\end{algorithm}
When $m_p \ll n_p $, the complexity of Algorithm~\ref{algo: tolerance} might be satisfactory.
Alternatively, complexity could be reduced to 
$O(n_p \log(m_p))$ by storing the $B_k^{-1}$ matrices in some sorted fashion (e.g., based on $||B_k||_2$).
A faster logarithmic search is then accomplished by first identifying a subset of $\database$ matrices to search and
then changing the {\sf for} loop on {\sf Line 5} so that only matrices in this subset are tested. 
Additionally, while this algorithm as written exits upon the first found match,
such behavior makes it sensitive to permutations of the input data in cases with large $\varepsilon$.
Instead, one may match $A_i$ with the closest of many $B_j$ satisfying {\sf Line 6}.
These modifications and other forms of acceleration for this approach will be considered in future work.

\subsection{Patch Clustering}
\label{subsec: clustering}

As noted, Algorithm~\ref{algo: tolerance} can be adapted so that it is relatively inexpensive.
However, Algorithm~\ref{algo: tolerance} prioritizes the minimization of accuracy term in Equation~\ref{eq: patch equality}.
When the number of database entries is a constraint, and the first term of $\mathcal{L}(\mathcal{B},\phi)$ is main priority for minimization,
instead one may consider clustering algorithms where the number of clusters is specified a priori.
In this paper, we consider perhaps one of the simplest based on $k$-means clustering.
$k$-means starts with an initial set of cluster assignments and computes an associated mean for each cluster.
The algorithm alternates between updating the cluster assignments
by identifying the points closest to the mean and then by updating the means given
the new cluster assignments.
In our context, we define the distance metric on patch matrices by 
\begin{align}
  d(A,B)
  &=
  \| I - AB^{-1} \|_2.
  \label{eq:clusterdistance}
\end{align}
Cluster means are then computed by the entry-wise mean
\begin{align}
  B_i
  &=
  |\phi^{-1}(i)|^{-1}
  \sum_{k\in\phi^{-1}(i)} A_k,
  \label{eq:clustermean}
\end{align}
where $\phi^{-1}(i)$ is the pre-image of $i$ under $\phi$ and $|\cdot|$ is set cardinality.
Because the map $\phi$ must be onto, this pre-image is always well-defined.
Algorithm~\ref{algo: clustering} summarizes the $k$-means algorithm.
\begin{algorithm}[H]
\caption{Spectral Clustering Construction of $\database$ and $\phi()$} \label{algo: clustering}
\begin{algorithmic}[1]
\State \textbf{Input:}   patches $\{A_1,A_2,\dots,A_{n_p}\}$, database size $m_p$
\State $\vec{\phi}=\mbox{randperm}(n_p,m_p)$;~~~~~\codeComment{Initialize clusters via random permutation}
\For{$j=1,\dots,m_p$}
  \State $B_j^{-1} = A_{\phi(j)}^{-1}$;~~~~~~~~~~~~~~\codeComment{Assign cluster representatives}
\EndFor
\State converged:=false;
\While{converged==false}
  \State converged=true;
  \For{$i=1,\dots,n_p$}~~~~~~~~~~~\codeComment{For each patch...}
    \For{$j=1,\dots,m_p$}~~~~~\codeComment{For each cluster...}
      \State $d_{j} = d(A_i,B_j)$;~~~~~~~~\codeComment{compute distance to cluster center}
    \EndFor
    \If{$\phi(i)\neq\mbox{argmin}_j(d_{j})$};~~~~\codeComment{If nearest cluster changed...}
      \State converged=false;~~~~~~~~~~~~~~\codeComment{clustering is not converged,}
      \State $\phi(i)=\mbox{argmin}_j(d_{j})$;~~~~~~~~~~\codeComment{update closest cluster assignment}
    \EndIf
  \EndFor
  \If{converged==false}~~~~~~~~\codeComment{If not converged...}
    \For{$j=1,\dots,m_p$}
      \State recompute $B_j^{-1}$;~~~~~~~~~~~~\codeComment{update cluster representatives}
    \EndFor
  \EndIf
\EndWhile
\State \textbf{Output:} $\mathcal{B}=\{B_1^{-1},B_2^{-1},\dots,B_{m_p}^{-1}\}$, $\vec{\phi}$
\end{algorithmic}
\end{algorithm}
We note that Algorithm~\ref{algo: clustering} provides one simple method for determining when to end the algorithm,
but there are other possible choices, such as maximum iteration limits, which are trivial to implement.
Line 11 defines the distance between patches, and line 20 refers to recomputing $B_i^{-1}$, which can be done several different ways.
We will break this algorithm into three sub-algorithms depending on how this is addressed:
\begin{itemize}
  \item \textbf{Entrywise $k$-means} computes $d(A_i,B_j) = \|A_i - B_j \|_{\ell_1}$ and defines $B_j^{-1}$ as the inverse of the entrywise average over a cluster.
  \item \textbf{Spectral $k$-means} computes $d(A_i,B_j) = \|I - A_i B_j^{-1} \|_2$ and defines $B_j^{-1}$ as the inverse of the entrywise average over a cluster.
  \item \textbf{Variance-minimizing clustering} computes $d(A_i,B_j) = \|I - A_i B_j^{-1} \|_2$ and defines $B_j^{-1}$ as the inverse of the cluster member which minimizes the in-cluster variance.
  \item \textbf{Bootstrapping} initializes clustering by first computing Algorithm~\ref{algo: tolerance} and then uses the resulting $\mathcal{B}$ as initial clusters for Algorithm~\ref{algo: clustering}. We expect this to perform slightly better than Algorithm~\ref{algo: tolerance} on its own.
\end{itemize}

As the algorithm is written, there are potential difficulties with computing cluster assignments in this fashion,
and future work will consider alternatives. One obvious problem
concerns the preservation of boundary conditions that sometimes appear in local patches.
Boundary conditions appear in these patches as a row with only one nonzero entry, and therefore any entrywise 
averaging scheme between a matrix with a boundary condition and a matrix without a boundary condition will immediately
result in a matrix without a boundary condition, destroying the original boundary condition information.
In practice, we adapt our method by partitioning the patch matrices dataset based on the presence of a
boundary condition, which splits the dataset into an ``interior'' and a ``boundary.''
Again, due to similar issues with entrywise averaging, we partition the boundary further so that
each boundary partition only contains matrices which have boundary conditions in the same rows.
In theory, there are possibly $2^{p_s}$ different partitions of the dataset;
however, in practice there are 9 on a two-dimensional square domain (associated with four corners, four edges, and one interior
of the domain)
and 27 on a three-dimensional cube domain (associated with eight corners, twelve edges, six faces, and one interior of the domain).
We then independently cluster on each partition,
and the final cluster assignments are simply the union of
all partitioned clustering results.
In Section~\ref{sec: numerical results}, we divide the number of clusters between boundary types
so that the number of patches is roughly proportional to the number of patches in the group.
One possible nearly equivalent alternative to the partitioning approach
would be to redefine the distance metric with an additional 
term that heavily penalizes cases where the number of nonzeros
in a row of $A$ and a row of $B$ are different.
This could be used in conjunction with some initialization of $\phi()$
to ensure that no initial clusters violate this nonzero condition. 

There are several possible improvements that one can consider to the basic $k$-means algorithms to
reduce the overall cost, automatically determine a cluster size, and to improve to the quality of the clusters 
found via the algorithm.
The Bootstrapping method effectively determines the number of clusters for $k$-means and greatly reduces the number of $k$-means iterations needed.
It also eliminates the need to partition patches based on their interior or boundary orientation.
Furthermore, cost may be greatly improved by performing minibatch $k$-means, although the algorithm
is typically less stable and depends on an adequate representation of the group by the random sample \cite{baraldi1999survey}.
In addition, there are a number of advanced possibilities that are somewhat more closely tied to our application.

For linear problem sequences (e.g., from
time marching or nonlinear solvers), there are techniques to update or revise an existing set of clusters
to account for the new data.
For example, one may initialize a database using the results from the previous time step
to greatly accelerate the convergence of the clustering algorithm.
However, this also depends on other factors such as the nonlinearity of the system.

Additionally, there are interesting combinations of clustering algorithms with 
deep learning neural networks that may be useful to compress the representation before the clustering algorithms,
determine the number of clusters, and allow for more general type of cluster regions to be determined~\cite{caron2018deep}.
We have initiated some work in the deep learning area, but this is beyond the scope of the current manuscript.

\section{Numerical Results}
\label{sec: numerical results}

\begin{figure}[H]
\centering
\resizebox{!}{5cm}{\includegraphics{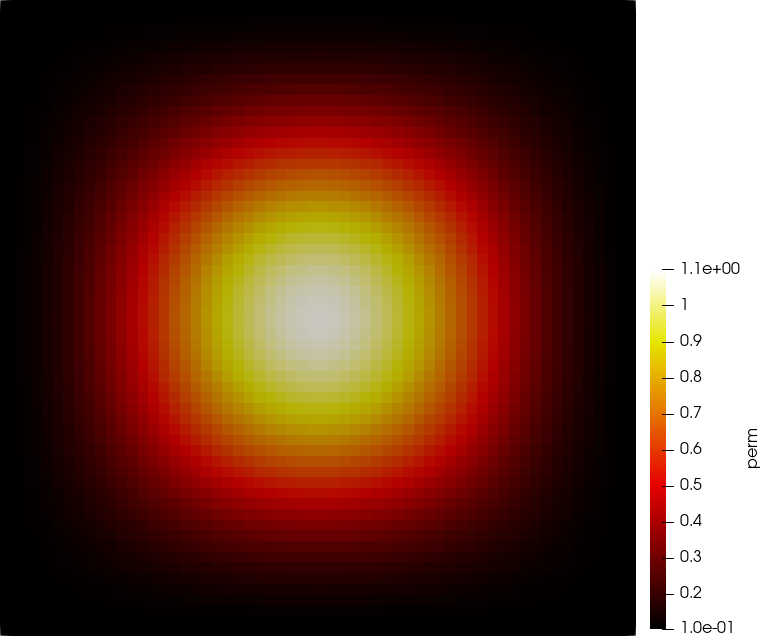}}
\resizebox{!}{5cm}{\includegraphics{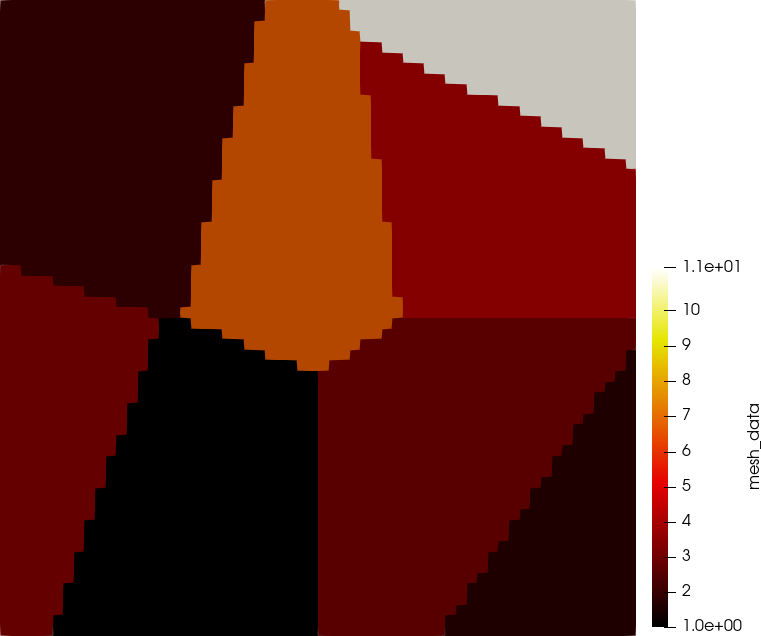}}
\caption{Graphical illustration of $\rho(x,y)$ for a smooth case (left) and for a discontinuous case (right)}
\label{fig: Poisson coefficients}
\end{figure}

We now consider several different scenarios to demonstrate our ability to identify and exploit structure in various scenarios.
For multigrid applications, a two dimensional Poisson problem is considered with varying coefficients
$$
-\nabla \cdot ( \rho(x,y) \nabla u )  = f
$$
on a unit square $[0,1]^2$ with Dirichlet boundary conditions.
Finite elements are employed on a regular quadrilateral mesh using
Lagrange polynomial basis functions with polynomial degree $p \le 5$.
Example 1 corresponds to a smoothly varying situation where $\rho()$ is given by a product of sine functions with a constant offset from zero and is depicted on the left side of Figure~\ref{fig: Poisson coefficients}.
Example 2, shown on the right side of Figure~\ref{fig: Poisson coefficients}, corresponds to a multiple material case where $\rho()$ is a piecewise constant function.
Example 3 corresponds to the case $\rho(x,y)=1$ and is used primarily for the performance discussion,
and Example 4 presents discussion on how one may consider the viability of extending such an algorithm to more complex problems such as Burgers' equation with a shock formation.

\textbf{Experiment 1} patch compression for the smooth case.
\\
Experiments 1 and 2 are performed in Matlab using GMRES (with a restart size of 20) in conjunction 
with $ M^{-1}_{combo} $. 
Again, each patch corresponds to one finite element,
so there are $3600$ true patches on a $60 \times 60$ mesh.
A direct solver applied to the matrix $R_0 A P_0$ is used for $M^{-1}_{amg}$,
though we have also performed experiments with AMG that give similar convergence profiles.
Here, $A$ is the high-order discretization matrix.
The interpolation operator $P_0$ corresponds to
linear interpolation viewing the spatial location of the DoFs 
in a high-order discretization as defining a fine mesh.
Finally, $R_0 = P_0^T$.
The function $\rho(x,y)$ is defined by $\rho(x,y) = \sin^2( \pi x ) \sin^2( \pi y ) + 0.1$.
The right hand side is chosen so that the solution is given by 
$ u = \sin(\pi x)\sin(\pi y)$. The GMRES iterative process is started with a zero initial guess
and convergence is declared when the initial residual is reduced by $10^{-8}$.
Table~\ref{table: example 1}
gives results for a variety of algorithms from Section~\ref{sec: metrics} applied to Example 1
on a $60 \times 60$ mesh.
and clustering is not run when the database size is less than 10
as the total number of clusters is split across the 8 boundary conditions and 1 interior,
as indicated by $-$ entries.

\begin{table}
  \centering
  \begin{tabular}{|l|l|l|l|l|l|l|l|l|l|}
  \hline
    \multirow{5}{*}{\rotatebox[origin=c]{90}{$p=2$}}
    & Algorithm\textbackslash Database Size & 3600 & 74 & 35 & 18 & 15 & 13 & 7 & 6 \\ \cline{2-10}
    & Greedy Tolerance~\ref{algo: tolerance} & 11 & 12 & 13 & 14 & 17 & 26 & 61 & 68 \\ \cline{2-10}
    & Spectral $k$-means ~\ref{algo: clustering} & 11 & 12 & 12 & 14 & 16 & 17 & $-$ & $-$ \\ \cline{2-10}
    & Var-Minimizing Clustering & 11 & 12 & 12 & 13 & 14 & 15 & $-$ & $-$ \\ \cline{2-10}
    & Entrywise $k$-means & 11 & 11 & 12 & 12 & 14 & 16 & $-$ & $-$ \\ \cline{2-10}
    \hline
    \hline
    \multirow{5}{*}{\rotatebox[origin=c]{90}{$p=3$}}
    & Algorithm\textbackslash Database Size & 3600 & 71 & 34 & 18 & 15 & 13 & 8 & 6 \\ \cline{2-10}
    & Greedy Tolerance~\ref{algo: tolerance} & 12 & 13 & 13 & 15 & 18 & 26 & 58 & 69 \\ \cline{2-10}
    & Spectral $k$-means ~\ref{algo: clustering} & 12 & 12 & 12 & 15 & 17 & 18 & $-$ & $-$ \\ \cline{2-10}
    & Var-Minimizing Clustering & 12 & 12 & 12 & 13 & 14 & 15 & $-$ & $-$ \\ \cline{2-10}
    & Entrywise $k$-means & 12 & 12 & 12 & 12 & 14 & 16 & $-$ & $-$ \\ \cline{2-10}
    \hline
    \hline
    \multirow{5}{*}{\rotatebox[origin=c]{90}{$p=4$}}
    & Algorithm\textbackslash Database Size & 3600 & 71 & 34 & 18 & 15 & 13 & 8 & 7 \\ \cline{2-10}
    & Greedy Tolerance~\ref{algo: tolerance} & 14 & 14 & 15 & 17 & 20 & 29 & 68 & 67 \\ \cline{2-10}
    & Spectral $k$-means ~\ref{algo: clustering} & 14 & 14 & 14 & 17 & 19 & 20 & $-$ & $-$ \\ \cline{2-10}
    & Var-Minimizing Clustering & 14 & 14 & 14 & 15 & 16 & 17 & $-$ & $-$ \\ \cline{2-10}
    & Entrywise $k$-means & 14 & 14 & 14 & 14 & 16 & 18 & $-$ & $-$ \\ \cline{2-10}
    \hline
    \hline
    \multirow{5}{*}{\rotatebox[origin=c]{90}{$p=5$}}
    & Algorithm\textbackslash Database Size & 3600 & 73 & 35 & 18 & 15 & 13 & 9 & 7 \\ \cline{2-10}
    & Greedy Tolerance~\ref{algo: tolerance} & 15 & 16 & 17 & 19 & 22 & 32 & 78 & 76 \\ \cline{2-10}
    & Spectral $k$-means ~\ref{algo: clustering} & 15 & 16 & 16 & 19 & 21 & 22 & $-$ & $-$ \\ \cline{2-10}
    & Var-Minimizing Clustering & 15 & 15 & 16 & 17 & 18 & 19 & $-$ & $-$ \\ \cline{2-10}
    & Entrywise $k$-means & 15 & 15 & 15 & 16 & 18 & 20 & $-$ & $-$ \\ \cline{2-10}
    \hline
    \hline
  \end{tabular}
  \caption{Comparison of multigrid iteration counts for varying database algorithms, database sizes, and polynomial degrees applied to Example 1 on a $60\times60$ mesh.}
  \label{table: example 1}
\end{table}

By examining Table~\ref{table: example 1}
one can see that with only 34 or 35 patches
the total number of GMRES iterations increases
by no more than 2 iterations over the case when all 3600 patches are used,
regardless of the choice of $p$.
That is, acceptable convergence rates are attained
even while retaining less than one percent of all patches.
However, as the database size is reduced to 10 patches or less,
convergence begins to degrade significantly for Algorithm~\ref{algo: tolerance}.

\textbf{Experiment 2} patch compression for the discontinuous case.
\\
The right hand side is chosen so that the solution is again given by 
$ u = \sin(\pi x)\sin(\pi y)$, 
$\rho$ is discontinuous piecewise constants,
and the iterative GMRES process is started with a zero initial guess and a tolerance of $10^{-8}$. 
Table~\ref{table: example 2}
gives multigrid iteration counts for a variety of algorithms from Section~\ref{sec: metrics}, database sizes,
and polynomial degrees applied to Experiment 2 on a $60 \times 60$ mesh,
and again clustering is not run when the database size is less than 10.

\begin{table}
  \centering
  \begin{tabular}{|l|l|l|l|l|l|l|l|l|l|}
  \hline
    \multirow{5}{*}{\rotatebox[origin=c]{90}{$p=2$}}
    & Algorithm\textbackslash Database Size & 3600 & 131 & 113 & 96 & 52 & 25 & 5 & 3 \\ \cline{2-10}
    & Greedy Tolerance~\ref{algo: tolerance} & 11 & 12 & 12 & 14 & 17 & 22 & 47 & 52 \\ \cline{2-10}
    & Spectral $k$-means ~\ref{algo: clustering} & 11 & 17 & 19 & 19 & 26 & 27 & $-$ & $-$ \\ \cline{2-10}
    & Var-Minimizing Clustering & 11 & 20 & 20 & 20 & 29 & 31 & $-$ & $-$ \\ \cline{2-10}
    & Entrywise $k$-means & 11 & 12 & 17 & 17 & 26 & 31 & $-$ & $-$ \\ \cline{2-10}
    \hline
    \hline
    \multirow{5}{*}{\rotatebox[origin=c]{90}{$p=3$}}
    & Algorithm\textbackslash Database Size & 3600 & 131 & 114 & 100 & 59 & 32 & 8 & 5 \\ \cline{2-10}
    & Greedy Tolerance~\ref{algo: tolerance} & 12 & 12 & 13 & 14 & 17 & 26 & 48 & 50 \\ \cline{2-10}
    & Spectral $k$-means ~\ref{algo: clustering} & 12 & 20 & 26 & 26 & 34 & 35 & $-$ & $-$ \\ \cline{2-10}
    & Var-Minimizing Clustering & 12 & 24 & 24 & 24 & 32 & 32 & $-$ & $-$ \\ \cline{2-10}
    & Entrywise $k$-means & 12 & 13 & 19 & 19 & 29 & 33 & $-$ & $-$ \\ \cline{2-10}
    \hline
    \hline
    \multirow{5}{*}{\rotatebox[origin=c]{90}{$p=4$}}
    & Algorithm\textbackslash Database Size & 3600 & 130 & 114 & 103 & 70 & 36 & 12 & 9 \\ \cline{2-10}
    & Greedy Tolerance~\ref{algo: tolerance} & 14 & 14 & 15 & 15 & 19 & 28 & 51 & 53 \\ \cline{2-10}
    & Spectral $k$-means ~\ref{algo: clustering} & 14 & 23 & 32 & 32 & 33 & 38 & 36 & $-$ \\ \cline{2-10}
    & Var-Minimizing Clustering & 14 & 27 & 27 & 27 & 27 & 38 & 37 & $-$ \\ \cline{2-10}
    & Entrywise $k$-means & 14 & 16 & 22 & 22 & 22 & 35 & 37 & $-$ \\ \cline{2-10}
    \hline
    \hline
    \multirow{6}{*}{\rotatebox[origin=c]{90}{$p=5$}}
    & Algorithm\textbackslash Database Size & 3600 & 130 & 116 & 103 & 78 & 38 & 15 & 10 \\ \cline{2-10}
    & Greedy Tolerance~\ref{algo: tolerance} & 15 & 16 & 16 & 17 & 21 & 29 & 60 & 64 \\ \cline{2-10}
    & Spectral $k$-means ~\ref{algo: clustering} & 15 & 25 & 35 & 30 & 35 & 40 & 43 & 39 \\ \cline{2-10}
    & Var-Minimizing Clustering & 15 & 28 & 28 & 54 & 47 & 40 & 40 & 44 \\ \cline{2-10}
    & Bootstrapped Var-Minimizing & 15 & 16 & 16 & 17 & 22 & 26 & 44 & 39 \\ \cline{2-10}
    & Entrywise $k$-means & 15 & 18 & 25 & 25 & 25 & 38 & 40 & 38 \\ \cline{2-10}
    \hline
    \hline
  \end{tabular}
  \caption{Comparison of multigrid iteration counts for varying database algorithms, database sizes, and polynomial degrees applied to Experiment 2 on a $60\times60$ mesh.}
  \label{table: example 2}
\end{table}

Interestingly, we see that Algorithm~\ref{algo: tolerance} outperforms any variation of clustering specified in Section~\ref{sec: metrics} except for the bootstrapped clustering, which we ran only in the case of $p=5$.
Bootstrapping performed best for $|\mathcal{B}| = 38, 15, 10$.
It required the same number of iterations as Algorithm~\ref{algo: tolerance} for other $|\mathcal{B}|$ sizes, except $|\mathcal{B}| = 78$ where it required one additional iteration.
With $|\database| = 130$ or 131, Algorithm~\ref{algo: tolerance} costs at most one additional solver iteration than using all 3600 patches,
while the clustering algorithms cost significantly more iterations and do not scale as well.
We suspect this is due to the piecewise constant nature of the problem,
which makes clustering less effective.

\begin{figure}[h]
\centering
\includegraphics[trim = 0.0in 0.0in 0.0in 0.0in, clip, height = 5.0cm,width = 5.0cm]{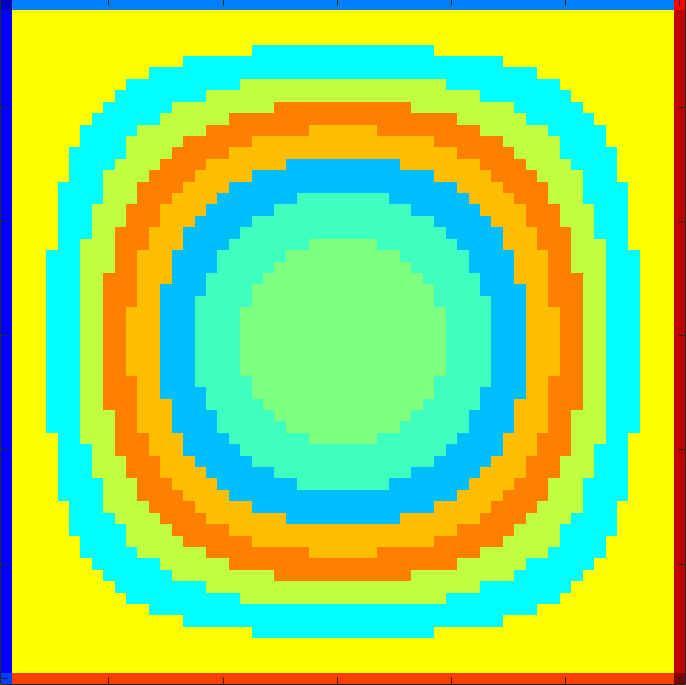}
~~~~
\includegraphics[trim = 0.0in 0.0in 0.0in 0.0in, clip, height = 5.0cm,width = 5.0cm]{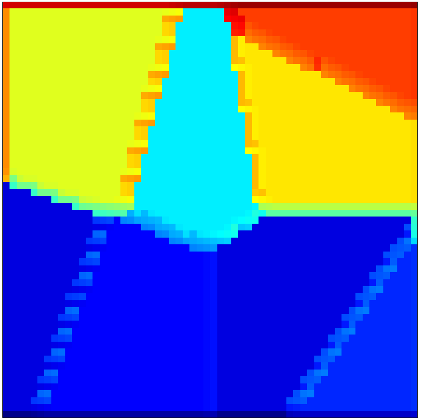}
\caption{Visualization of $\phi$ mapping for $p=5$ determined by Algorithm~\ref{algo: clustering} spectral $k$-means for Experiment 1 (left with 16 clusters) and by Algorithm~\ref{algo: tolerance} for Experiment 2 (right with 130 unique entries)}
\label{fig: pretty clusters}
\end{figure}

Figure~\ref{fig: pretty clusters} depicts the clusters identified by the two respective algorithms. Here, one can
see that the found clusters generally match the $\rho()$ functions depicted in 
Figure~\ref{fig: Poisson coefficients}.

\textbf{Experiment 3} Timings for $\rho(x,y)=1$ problem.
\\
While our focus is on reduction in the size of the dataset with little loss of accuracy, 
we did investigate timings using a variant of the greedy database compression scheme described by Algorithm~\ref{algo: tolerance}
in the \texttt{DatabaseSchwarz} class of the Ifpack2 preconditioning package, and we called it from the MueLu multigrid package.
Both packages are part of the C++ Trilinos framework.
We utilized this preconditioner as a simple smoothing iteration of the form
\begin{align*}
  x \leftarrow x + \omega M^{-1} r,
\end{align*}
where $r=b-Ax$ is the residual. 

For this, we use a variation of Algorithm~\ref{algo: tolerance} which checks $\|A_i - B_j \|_{\ell_1}<\varepsilon$ for simplicity.
In this case, we consider a $20\times20\times20$ finite element discretization with $p=5$ polynomial degree, yielding 8000 patches where each patch is of size $216\times216$.
We utilize this approach as a preconditioner inside a GMRES solver with a database tolerance of $\varepsilon=10^{-7}$ and a GMRES tolerance of $10^{-7}$ as well.
We compare timings in two different scenarios, each taking 39 iterations to solve.
First, no compression is applied. In the no compression case, the setup phase computes LU factorizations for all 8000 patches, and the preconditioning apply phase solves the inversion for all 8000 patches using the previously computed LU factors.
Second, compression is applied. In this case, the setup phase starts by using Algorithm~\ref{algo: tolerance} to construct $\mathcal{B}$ with the aforementioned modification.
For this problem, $\mathcal{B}$ has exactly 27 entries.
After the database is computed, the LU factors are stored in-place instead of strict inverses.
The apply phase then applies the 27 stored factorizations against the appropriate choices for the 8000 patches.
In each case, we use the LAPACK GETRF and GETRS routines to compute and apply these factorizations, respectively.
Table~\ref{table: timings} compares timings for the two approaches on a node of the Attaway supercomputer at Sandia.

\begin{table}
  \centering
  \begin{tabular}{|l|l|l|l|l|}
  \hline
    Configuration & Setup time (s) & Apply time (s) \\
  \hline
  No compression & 287.4 & 64.5 \\
  \hline
  Compression & 262.0 & 49.4 \\
  \hline
  \end{tabular}
  \caption{Experiment 3: comparison of timings in the case of no compression and compression}
  \label{table: timings}
\end{table}

We observe a slight (9\%) speedup in the setup phase of the {\sf Compression} case because it computes 27 factorizations as opposed to 8000 factorizations,
which outweighs the cost of repeatedly computing $\|A_i - B_j \|_{\ell_1}$ for the compression.
We suspect that our current implementation of Algorithm~\ref{algo: tolerance} could be improved for larger database sizes by using a logarithmic search
algorithm as discussed earlier. Despite these limitations, there is a significant (23.4\%) speedup in the 
apply phase when compression is utilized even with our simple implementation. 
We believe this is due to increased data locality by accessing only 27 factorizations of size $216\times216$ as opposed to accessing 8000 factorizations of size $216\times216$.
Assuming 8 bytes per value,
this corresponds to a difference of storing approximately 9.6MB of dense patch factorizations and 62.5KB of database lookup indices in the compressed case versus 2.8GB of dense patch factorizations in the uncompressed case.
To put this in perspective, the cheapest smoother to store explicitly is the Jacobi method, which would require storing 7.9MB of data to store the diagonal of the $1030301\times 1030301$ matrix for this problem.
We also remark that our implementation has plenty of room for performance enhancements,
and these numbers should be considered a lower bound for potential benefits. For example, we believe it would be much better to 
simultaneously solve all patch problems within the same cluster (which we are not currently doing) and the patch matrix could
be loaded once in memory and then used for the multiple right hand sides. 
Additionally, we expect problems with larger patches will be solved with even
larger performance gains over approaches with no compression due to the cost scaling of inverting patches.
Finally, we note that for linear problem sequences (e.g., within nonlinear solvers or time stepping schemes), setup
times can probably be amortized by incorporating some reuse of database assignments from the previous linear setup.

\begin{figure}[h]
\centering
\resizebox{!}{3.5cm}{
\includegraphics{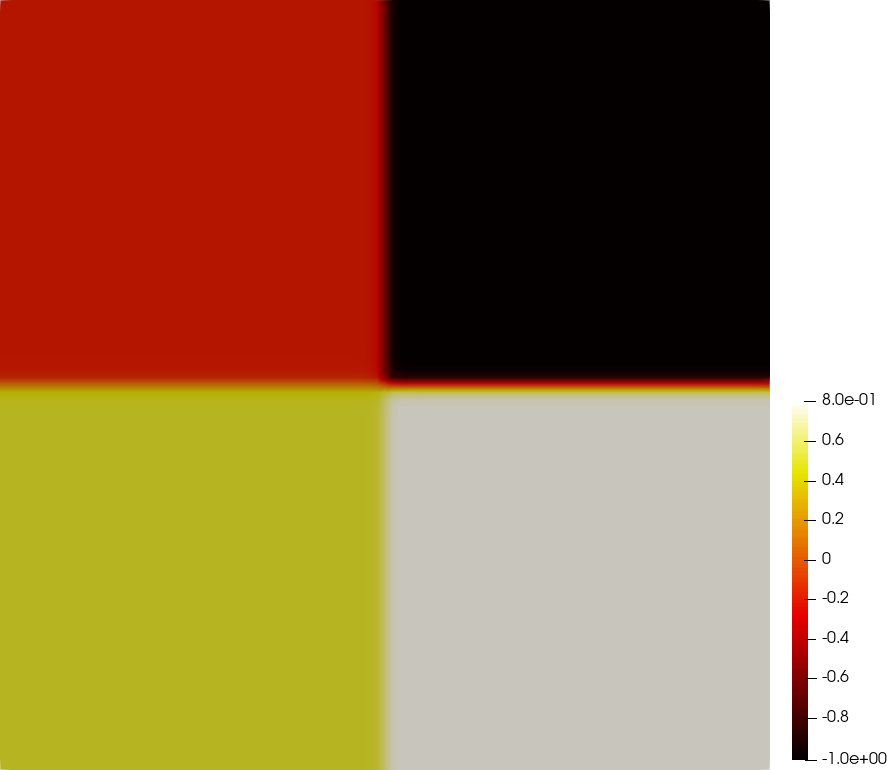}}
\resizebox{!}{3.5cm}{
\includegraphics{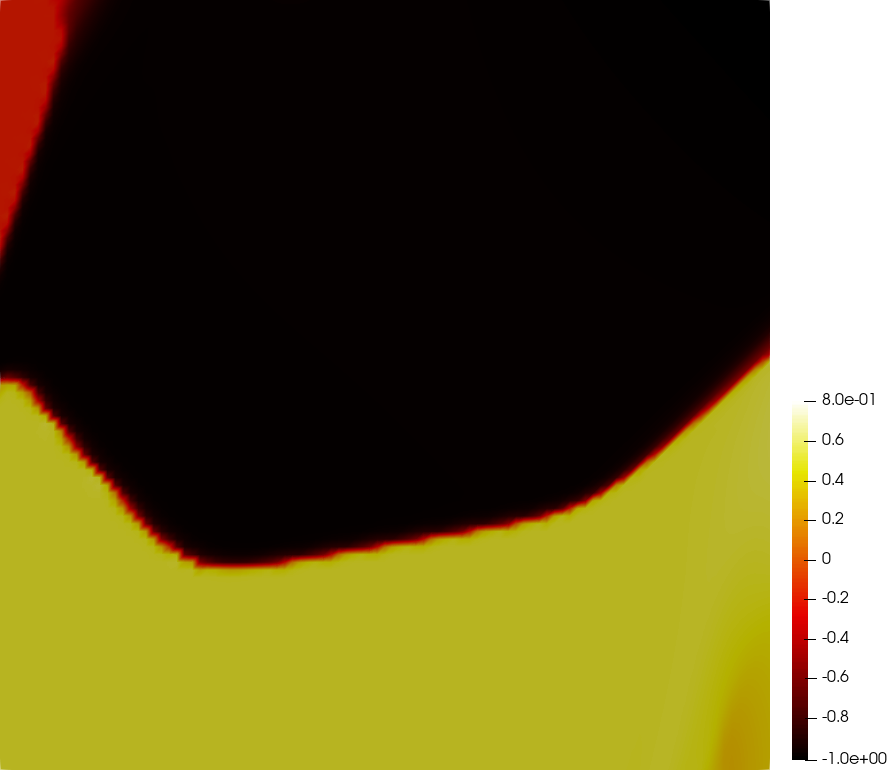}}
\resizebox{!}{3.5cm}{
\includegraphics{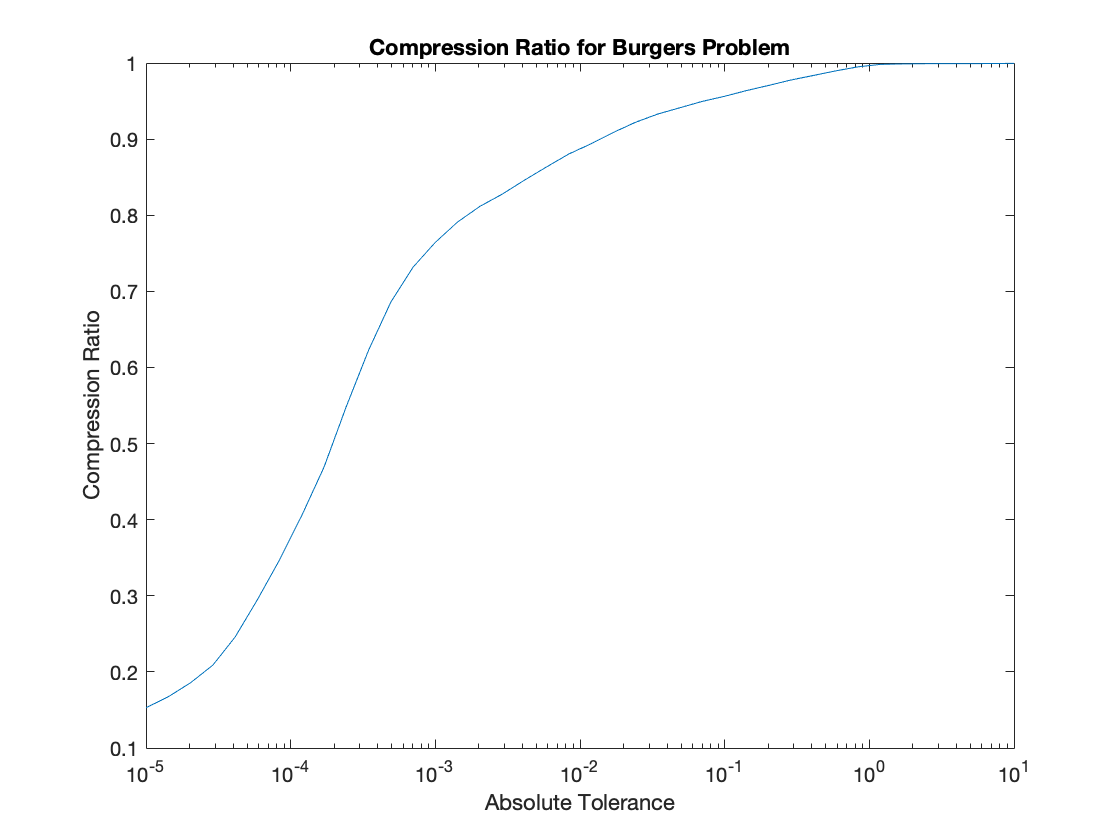}}
\caption{Burgers' equation from initial condition $t=0$ (left) to shock formation at final time $t=1$ (center) and the associated compressibility curve at final time based on Algorithm~\ref{algo: tolerance} (right) on a $100\times100$ mesh with $p=2$.}
\label{fig: burgers}
\end{figure}

\textbf{Experiment 4} Nonlinear shock dynamics.
We close this section with an example of a nonlinear system with a shock formation.
This example utilizes Burgers' equation
\begin{align*}
  \frac{du}{dt} + \nabla\cdot \left(\frac{1}{2} \vec{\nu} u^2 - \epsilon(\vec{x}) \nabla u \right) = f(\vec{x}).
\end{align*}
with advection $\vec{\nu}$, entropy viscosity strength $\epsilon$.  Additional terms for SUPG stabilization are present but not shown.
We focus on the case where a shock develops in the primary variable $u$, using $\vec{\nu}=[1,1]^T$. 
At the final time indicated by the center image of Figure~\ref{fig: burgers},
we extract the matrix and analyze the system for compressibility using Algorithm~\ref{algo: tolerance}.

In particular, we do not solve the system;
however, we use this experiment to demonstrate how patch compressibility behaves for much more complex nonlinear physics.
That is, the right image is a plot of $\frac{n_p-\|\mathcal{B}\|}{n_p}$ for a wide range of values of $\varepsilon$.
It is important to note that compression ratios as high as 90\% may be achieved
with a tolerance on the order of $\varepsilon=10^{-2}$,
meaning this approach is likely applicable to a broad range of problems.

This approach may be extended to Algorithm~\ref{algo: clustering}
by plotting a clustering quality measure such as cluster variance as a function of number of clusters,
which will help identify how amenable a system is to the algorithms in this paper.
Additionally not shown, one may plot a histogram of the mapping $\phi()$ to identify the distribution of structure for a given matrix.
For example, a histogram for $\phi()$ related to Experiment 2 would show large peaks corresponding to the interior of the constant regions,
indicating the problem is largely structured with the exception of interfaces.

\section{Conclusion}
\label{sec: conclusion}
Two families of algorithms have been presented for detecting and exploiting structure in a linear system for patch preconditioners based on the idea that similar patches can share the
same factorization when solving the small patch subproblems. We have illustrated that on some examples, it is possible to maintain
similar convergence rates even when the number of stored/factored patches is less than five percent of the number of true
patches. 
This gain makes a more accurate patch-based approach more competitive with
inexpensive Jacobi methods, which require very little additional storage (specifically only the matrix diagonal).
This is essential in matrix-free applications where often a matrix-free approach is adopted due significant
storage concerns.
Previously, matrix-free storage advantages are often lost when all patches are stored/factored,
as the memory needed to store all patch factorizations is
often comparable to the storage needed for the discretization matrix.
Now, by only factoring a small subset of the patch matrices,
the storage advantages of a matrix-free approach are once again possible.
While our focus has been on storage reduction,
we have presented some results hinting that it may be possible to reduce
run time using a patch compression approach even when a matrix-free approach is not being used.
However, additional work is needed to truly demonstrate this.
The algorithms presented in this paper are limited
in terms of both their implementation and their sophistication.
More sophisticated algorithms should be considered
to reduce the run time during the setup phase when defining
the patch clusters and to reduce the run time
in the apply phase when repeatedly solving the same patch matrix system with
many different right hand sides.
Additionally, more sophisticated clustering algorithms should be considered to improve 
the quality of the clustering perhaps leveraging ideas
that combine clustering algorithms with deep learning algorithms.

The authors would like to thank Nicholas Moore for valuable discussions and feedback.
This work was supported by the U.S.~Department of Energy, Office of Science, Office of Advanced Scientific
Computing Research, Applied Mathematics program and Sandia National Laboratories Laboratory Directed
Research \& Development Program.
Sandia National Laboratories is a multimission laboratory managed and operated by National Technology
and Engineering Solutions of Sandia, LLC., a wholly owned subsidiary of Honeywell International, Inc.,
for the U.S. Department of Energy's National Nuclear Security Administration under grant~DE-NA-0003525.
This paper describes objective technical results and analysis.  Any subjective views or opinions that
might be expressed in the paper do not necessarily represent the views of the U.S. Department of Energy
or the United States Government.

\bibliographystyle{siam}
\bibliography{patch_smoothers.bib}

\end{document}